\documentclass[12pt]{article}
\usepackage{amssymb}
\usepackage{latexsym,bm}
\usepackage{graphicx}
\usepackage{amsmath}
\usepackage{mathrsfs}
\usepackage{epstopdf}

\date{}
\newcounter{mathitem}
  {\begin{list}{{$(\roman{mathitem})$}}{
   \setcounter{mathitem}{0}
   \usecounter{mathitem}
   \setlength{\topsep}{0pt plus 2pt minus 0pt}
   \setlength{\parskip}{0pt plus 2pt minus 0pt}
   \setlength{\partopsep}{0pt plus 2pt minus 0pt}
   \setlength{\parsep}{0pt plus 2pt minus 0pt}
   \setlength{\leftmargin}{35pt}
   \setlength{\itemsep}{0pt plus 2pt minus 0pt}}}
  {\end{list}}

\begin{document}
\title{Cycles of consecutive lengths in $3$-connected graphs
\footnote{E-mail addresses: {\tt lichengli0130@126.com} (C. Li), {\tt zhan@math.ecnu.edu.cn} (X. Zhan).}}
\author{\hskip -10mm Chengli Li and Xingzhi Zhan\thanks{Corresponding author}\\
{\hskip -10mm \small Department of Mathematics,  Key Laboratory of MEA (Ministry of Education) }\\
{\hskip -10mm \small \& Shanghai Key Laboratory of PMMP, East China Normal University, Shanghai 200241, China}}\maketitle
\begin{abstract}
Recently Lin, Wang and Zhou have proved that every $3$-connected nonbipartite graph of minimum degree at least $k$ with $k\ge 6$ and order at least
$k+2$ contains $k$ cycles of consecutive lengths. They also conjecture that this result is true for $k=4, 5.$ We prove this conjecture. Our proofs use many ideas
of Gao, Huo, Liu and Ma.
\end{abstract}


\vskip 8mm

\section{Introduction}

In 1998, Bondy and Vince [2] proved that every $3$-connected nonbipartite graph contains two cycles of consecutive lengths, and they
posed the following

{\bf Problem.} Does there exist a function $f(k)$ such that every $3$-connected nonbipartite graph with minimum degree at least $f(k)$ contains $k$ cycles of consecutive
lengths?

In 2002 Fan [4] solved this problem positively and proved that $f(k)$ may be taken to be $3\lceil k/2\rceil.$ Twenty years later Gao, Huo, Liu and Ma [5] proved
the following result which shows that $f(k)$ may be taken to be $k+1.$

{\bf Theorem 1} (Gao-Huo-Liu-Ma [5]). {\it Every $3$-connected nonbipartite graph with minimum degree at least $k+1$ contains $k$ cycles of consecutive lengths.}

Clearly, the nonbipartiteness condition is necessary in Theorem 1, since bipartite graphs have no odd cycles. Also, Bondy and Vince [2, p.12] constructed examples to show
that the $3$-connectedness condition is necessary.

Recently Lin, Wang and Zhou [8] have strengthened Theorem 1 by proving the following result. The {\it order} of a graph is its number
of vertices, and the {\it size} is its number of edges.

{\bf Theorem 2} (Lin-Wang-Zhou [8]). {\it Every $3$-connected nonbipartite graph of minimum degree at least $k$ with $k\ge 6$ and order at least $k+2$
contains $k$ cycles of consecutive lengths.}

They [8,  Conjecture 5.1] conjecture that Theorem 2 is also true for $k=4, 5.$ In this paper we prove this conjecture. Our proofs use many ideas
of Gao, Huo, Liu and Ma [5].

There are other interesting recent work on cycle lengths (e.g. [3] and [7]).

We will use the following standard notations. We denote by $V(G)$ and $E(G)$ the vertex set and edge set of a graph $G,$ respectively, and denote by $|G|$ and $||G||$ the order and size of $G,$ respectively. Thus, if $H$ is a cycle or a path then $||H||$ means the length of $H.$ The neighborhood and degree of a vertex $x$ in a graph $G$ is denoted by $N_G(x)$ and ${\rm deg}_G(x),$ respectively. If the graph $G$ is clear from the context, the subscript $G$ might be omitted.  We denote by $\delta(G)$ the minimum degree of $G.$ For a vertex subset $S\subseteq V(G),$ we use $G[S]$ to denote the subgraph of $G$ induced by $S.$ For two graphs $G$ and $H,$ $G\vee H$ denotes the {\it join} of $G$ and $H,$ which is obtained from the disjoint union $G+H$ by adding edges joining every vertex of $G$ to every vertex of $H.$

An $(x,y)$-path is a path with endpoints $x$ and $y.$ Suppose that $H$ is a subgraph or a vertex
subset of a graph $G.$ The neighborhood and closed neighborhood of $H$ in $G,$ denoted by $N_G(H)$ and $N_G[H],$ respectively are defined by
$$
N_G(H)=\left(\mathop{\cup}\limits_{x\in V(H)}N_G(x)\right)\setminus V(H),\quad N_G[H]=N_G(H)\cup V(H).
$$
A vertex $v$ and $H$  are  said to be {\it adjacent} if $v$ is adjacent to a vertex in $H.$

Sometimes we need to specify two special vertices in a graph. If $x$ and $y$ are two distinct vertices of a graph $G,$ we say that the triple
$(G,x,y)$ is a {\it rooted graph} (with roots $x$ and $y$). The {\it minimum degree} of a  rooted graph $(G,x,y)$ is defined to be ${\rm min}\{{\rm deg}_G(v)|\,v\in V(G)\setminus\{x,y\}\}.$ $(G,x,y)$ is said to be {\it $2$-connected} if $G+xy$ is $2$-connected, where $G+xy$ means $G$ if $x$ and $y$ are adjacent and $G+xy$ means the graph obtained from $G$ by adding the edge $xy$ if $x$ and $y$ are nonadjacent.

Let $S\subseteq V(G).$ The graph obtained from $G$ by {\it contracting $S$ to a vertex $s$} is the graph with vertex set $(V(G)\setminus S)\cup \{s\}$ and edge set
$E(G-S)\cup \{su|\,u\in N_G(S)\}.$

Recall that every connected graph $G$ has a block-cutvertex tree ([1, p.121], [10, p.156]), a leaf of which is called an {\it end-block}. If $G$ has cut-vertices, then
an end-block of $G$ contains exactly one cut-vertex.
If $B$ is an end-block and a vertex $b$ is the only cut-vertex of $G$ with $b\in V(B),$ then we say that $B$ is an {\it end block with cut-vertex $b$}. A {\it $k$-cycle}, denoted by $C_k,$ is a cycle of length $k.$ We denote by $K_s$ the complete graph of order $s.$
Let $K_{s,\, t}$ denote the complete bipartite graph whose partite sets have cardinality $s$ and $t,$ respectively. We denote by $\overline{G}$ the complement of a graph $G.$

\section{Main results}

The main result is as follows.

{\bf Theorem 3}. {\it Every $3$-connected nonbipartite graph of minimum degree at least $k$ with $k\ge 4$ and order at least $k+2$
contains $k$ cycles of consecutive lengths.}

We remark that the conclusion of Theorem 3 does not hold for $k=3.$ The Petersen graph is a $3$-connected nonbipartite graph of minimum degree $3$
whose cycle spectrum is $\{5, 6, 8, 9\}.$ Thus the Petersen graph does not contain three cycles of consecutive lengths.

To prove Theorem 3, we will need the following five lemmas. A set of paths in a graph is called {\it nice} if their lengths form an arithmetic progression with
initial term at least $2$ and with common difference $2.$

{\bf Lemma 4}. {\it Let $x$ and $y$ be two distinct vertices in a graph $G.$ If $G-x$ is triangle-free and $(G,x,y)$ is a $2$-connected rooted graph with minimum degree
at least $3,$ then $G$ contains two nice $(x,y)$-paths.}

{\bf Proof.} We use induction on the order of $G.$ The smallest order of such a graph is $5,$  and $\overline{K_2}\vee (K_2+K_1)$ and $K_1\vee C_4$
are the only two graphs of order $5$ that satisfy the conditions of Lemma 4. It is easy to verify that the conclusion of Lemma 4 holds for these two graphs. Next
let $G$ have order at least $6$ and assume that Lemma 4 holds for all graphs of orders less than $|G|.$

We may suppose $x$ and $y$ are nonadjacent in $G.$ Otherwise we may delete the edge $xy.$ We distinguish two cases.

Case 1. There exists a $4$-cycle containing $x$ in $G-y.$

Let $C=xx_1ax_2x$ be a $4$-cycle in $G-y,$ and let $F$ be the component of $G-V(C)$ that contains $y.$

Subcase 1.1. One of $x_1$ and $x_2,$ say $x_1,$ has a neighbor $z$ in $F.$

Let $P$ be a $(y,z)$-path in $F.$ Then the two paths $xx_1z\cup P$ and  $xx_2ax_1z\cup P$ are nice.

Subcase 1.2. None of $x_1$ and $x_2$ has a neighbor in $F.$

Since $(G,x,y)$ is $2$-connected, $a$ must have a neighbor $b$ in $F.$ Let $Q$ be a $(b,y)$-path in $F.$
Now $G-V(F)-x$ is triangle-free and $(G-V(F),x,a)$ is a $2$-connected rooted graph with minimum degree at least $3$ whose order is less than that of $G.$
By the induction hypothesis, $G-V(F)$ contains two nice $(x, a)$-paths $P_1$ and $P_2.$ Then  $P_1\cup ab\cup Q$ and $P_2\cup ab\cup Q$ are two nice $(x,y)$-paths in $G.$

Case 2. There exists no $4$-cycle containing $x$ in $G-y.$

Let $X=N_G(x).$ Note that by our assumption, $y\notin X.$ Let $G^*$ be the graph obtained from $G-x$ by contracting $X$ into a new vertex $x^*.$  Since $G-y$ does not have a $4$-cycle
containing $x,$ we deduce that for any vertex $v\in V(G^*)\setminus \{x^*, y\},$ $v$ has at most one neighbor in $X.$ It follows that ${\rm deg}_{G^*}(v)={\rm deg}_{G}(v).$

Denote $H=G^*+x^*y.$ If $H$ is not $2$-connected, then $x^*$ is the only cut-vertex of $H$ and every block of $H$ contains $x^*.$ In this case we let $B$ be the block of $H$
containing $y.$ If $H$ is $2$-connected, we let $B=H.$

If $|B|\ge 3,$ then $B$ is $2$-connected. Now $B-x^*$ is triangle-free and $(B, x^*, y)$ is a $2$-connected rooted graph with minimum degree at least $3.$
By the induction hypothesis, $B$ contains two nice $(x^*, y)$-paths, which yield $(u_i, y)$-path $P_i$ in $G$ with $u_i\in X$ for $i=1,2$ such that $P_1$ and $P_2$ are two nice paths where it is possible that $u_1=u_2.$ Then $xu_1\cup P_1$ and $xu_2\cup P_2$ are two nice $(x,y)$-paths.

It remains to consider the case when $|B|=2;$ i.e. $B=x^*y.$ Then $N_G(y)\subseteq X.$

Subcase 2.1. There exists a vertex $y^{\prime}\in N_G(y)$ such that $y^{\prime}$ is adjacent to $G-(X\cup\{x,y\}).$

Assume that $y'$ is adjacent to a component $D$ of $G-(X\cup\{x,y\}).$ Since $G+xy$ is $2$-connected, we have $|N_G(D)\cap X|\ge 2.$ Let $G_1$ be the graph obtained
from $G[X\cup D]$ by contracting $X\setminus\{y'\}$ into a new vertex $w_1.$ For any vertex $v\in V(D)$ we have ${\rm deg}_{G_1}(v)={\rm deg}_{G}(v).$
Moreover, $(G_1, w_1, y')$ is a $2$-connected rooted graph with minimum degree at least $3,$ and $G_1-w_1$ is triangle-free. By the induction hypothesis,
$G_1$ contains two nice $(w_1, y')$-paths, which yield $(u_i, y')$-path $P_i$ in $G$ with $u_i\in X\setminus \{y'\}$ for $i=1,2$ such that $P_1$ and $P_2$ are two nice paths where it is possible that $u_1=u_2.$ Then $xu_1\cup P_1\cup y'y$ and $xu_2\cup P_2\cup y'y$ are two nice $(x,y)$-paths.

Subcase 2.2. $N_G(y')\subseteq X\cup\{x,y\}$ for any vertex $y^{\prime}\in N_G(y).$

Recall that $G-x$ is triangle-free and we have assumed $x\notin N_G(y).$ Hence  $N_G(y)$ is an independent set. Note also that ${\rm deg}_G(y')\ge 3$ for any
$y^{\prime}\in N_G(y).$ Choose vertices $y^*\in N_G(y)$ and $y_1\in N_G(y^*)\setminus\{y,x\}.$ Then $y_1\in X$ and $y_1$ is nonadjacent to $y.$

If $y_1$ has a neighbor in $X$ other than $y^*,$ say, $y_2,$ then $xy^*y$ and $xy_2y_1y^*y$ are two nice $(x,y)$-paths. Otherwise $N_G(y_1)\cap X=\{y^*\},$ and next we make
this assumption.

Now $y_1$ is adjacent to $G-(X\cup\{x,y\}).$ Assume that $y_1$ is adjacent to a component $D_1$ of $G-(X\cup\{x,y\}).$ Since $G+xy$ is $2$-connected, we have $|N_G(D_1)\cap X|\ge 2.$ Let $G_2$ be the graph obtained from $G[X\cup D_1]$ by contracting $X\setminus \{y_1\}$ into a new vertex $w_2.$ For any vertex $v\in V(D_1)$ we have
${\rm deg}_{G_2}(v)={\rm deg}_{G}(v).$ Moreover, $(G_2, w_2, y_1)$ is a $2$-connected rooted graph with minimum degree at least $3,$ and $G_2-w_2$ is triangle-free. By the induction hypothesis, $G_2$ contains two nice $(w_2, y_1)$-paths, which yield $(u_i, y_1)$-path $P_i$ in $G$ with $u_i\in X\setminus \{y_1\}$ for $i=1,2$ such that $P_1$ and $P_2$ are two nice paths where it is possible that $u_1=u_2.$ Observe that $y^*\notin N_G(D_1).$ Thus $u_1\neq y^*$ and $u_2\neq y^*,$ and so  $xu_1\cup P_1\cup y_1y^*y$ and
$xu_2\cup P_2\cup y_1y^*y$ are two nice $(x,y)$-paths. \hfill$\Box$

We say that three paths $P_1, P_2, P_3$ are {\it good} if $||P_1||>||P_2||>||P_3||\ge 2$ and one of the following two conditions holds:
\newline\centerline{(1) $P_1,\, P_2,\, P_3$ are nice; \,\,\, (2) $\{||P_1||-||P_2||,\, ||P_2||-||P_3||\}=\{1, 2\}.$}

{\bf Lemma 5.} {\it
Let $G$ be a graph with $x,y\in V(G).$
If $(G,x,y)$ is a $2$-connected rooted graph with minimum degree at least $4$ and $G-x$ is triangle-free, then $G$ contains three good $(x,y)$-paths.
}

{\bf Proof.}
We use induction on the order of $G.$ The smallest order of such a graph is $7.$  $R=K_1\vee K_{3,3}$ and the graph obtained from $R$ by deleting one edge incident to the vertex in $K_1$ are the only two graphs of order $7$ that satisfy the conditions of Lemma 5. It is easy to verify that the conclusion of Lemma 5 holds for these two graphs. Next let $G$
have order at least $8$ and assume that Lemma 5 holds for all graphs of orders less than $|G|.$

We will use proof by contradiction. To the contrary, suppose that $G$ does not contain three good $(x,y)$-paths.

We may suppose $x$ and $y$ are nonadjacent in $G.$ Otherwise we may delete the edge $xy.$

{\bf Claim 1.}
$G$ is $2$-connected.

Since $G+xy$ is $2$-connected, $G$ is connected.
Suppose that $G$ is not $2$-connected.
Then $G$ contains a cut-vertex $b$ and two connected subgraphs $G_1,G_2$ of order at least two such that $G=G_1\cup G_2$ and $V(G_1)\cap V(G_2)=\{b\}$.
Note that exactly one of $x$ and $y$ lies in $G_1$ and the other lies in $G_2.$ Also $b\notin\{x,y\}.$ Without loss of generality we assume that $x\in V(G_1)\setminus\{b\}$,
 $y\in V(G_2)\setminus \{b\}.$ Clearly, either $|G_1|\geq 3,$ or $|G_2|\geq 3.$

Suppose that $|G_1|\geq 3.$ Then $(G_1,x,b)$ is a $2$-connected rooted graph with minimum degree at least $4$ and $G_1-x$ is triangle-free.
By the induction hypothesis, $G_1$ contains three good $(x,b)$-paths and so
concatenating each of these paths with a fixed $(b,y)$-path in $G_2$, we obtain three good $(x,y)$-paths in $G$, a contradiction.

Suppose that $|G_2|\geq 3.$ Then $(G_2,y,b)$ is a $2$-connected rooted graph with minimum degree at least $4$ and $G_2-y$ is triangle-free.
By the induction hypothesis, $G_2$ contains three good $(y,b)$-paths and so
concatenating each of these paths with a fixed $(b,x)$-path in $G_2$, we obtain three good $(x,y)$-paths in $G$, a contradiction.
Therefore $G$ is $2$-connected.
The proof of Claim~1 is complete.

{\bf Claim 2.} $G-y$ is triangle-free.

To the contrary, suppose that $G-y$ contains a triangle $T$. Since $G-x$ is triangle-free, $T$ contains the vertex $x$. Assume that $T=xx_1x_2.$
Let $C$ be the component of $G-V(T)$ containing $y.$

We assert $|C|\geq 2.$ Otherwise $V(C)=\{y\}$. Since $x$ and $y$ are nonadjacent in $G$ and $G$ is $2$-connected by Claim 1, we deduce that $y$ is adjacent to $x_1$ and $x_2$, which contradicts the fact that $G - x$ is triangle-free.

If $C$ is $2$-connected, then let $B=C$ and $b=y$;
otherwise let $B$ be an end-block of $C$ with cut-vertex $b$ such that $y\notin V(B)\setminus\{b\}$.

If $|B|=2,$ then $B$ has a vertex with degree one other than $y,$ say $b'.$ Since
$(G,x,y)$ has minimum degree at least 4, $b'$ is adjacent to each vertex in $T$ and so $b'x_1x_2$ is a triangle, a contradiction. Thus $|B|\ge 3$ and $B$ is $2$-connected. 

Since $G$ is $2$-connected by Claim~1, we have $N_G(B-b)\cap V(T)\neq \emptyset.$
If $N_G(B-b)\cap V(T)=\{x\},$ then let $G_1=G[V(B)\cup \{b,x\}].$ We observe that $(G_1,x,b)$ is a $2$-connected rooted graph with minimum degree at least $4$ and $G_1-x$ is triangle free. By the induction hypothesis, $G_1$ contains three good $(x,b)$-paths. Let $P$ be a $(b, y)$-path in $C-V(B-b).$ By concatenating these three paths with $P,$ we obtain three good $(x,y)$-paths in $G,$ a contradiction.

Thus we have $N_G(B-b)\cap \{x_1,x_2\}\neq \emptyset.$
Let $G_2$ be the graph obtained from $G[V(B)\cup \{x_1,x_2\}]$ by contracting $\{x_1,x_2\}$ into a vertex $x'$.
Recall that $G-x$ is triangle-free. Then for any vertex $v\in V(B-b),$ we have $\deg_{G_2}(v)\ge \deg_{G}(v)-1\ge 3.$
Clearly, $(G_2,x',b)$ is a $2$-connected graph with minimum degree at least $3$ and $G_2-x'$ is triangle-free. By Lemma 4, $G_2$ contains two nice $(x',b)$-paths, which yield $(u_i,b)$-path $P_i$ in $G$ with $u_i\in \{x_1,x_2\}$ for $i=1,2$ such that $||P_1||=||P_2||+2,$ where it is possible that $u_1=u_2.$ So $xu_2P_2bPy,$ $xu_1P_1bPy$ and $xu_2u_1P_1bPy$ are three good $(x,y)$-paths in $G,$ a contradiction.
The proof of Claim~2 is complete.

Recall that $x$ and $y$ are not adjacent. By Claim~2,  $G$ is triangle-free.
In what follows, due to the symmetry of $x$ and $y$ in $G$, we can assume without loss of generality that ${\rm deg}_{G}(x)\leq \deg_G(y)$.

{\bf Claim 3.}
There exists a $4$-cycle containing $x$ in $G-y.$

To the contrary, suppose that there exists no $4$-cycle containing $x$ in $G-y$. Let $X=N_G(x).$ Note that by our assumption, $y\notin X.$
Let $G_1$ be the graph obtained from $G-x$ by contracting $X$ into a new vertex $x_1$. Since $G-y$ does not have a $4$-cycle containing $x,$ we deduce that for any vertex $v\in V(G_1)\setminus \{x_1,y\},$ $v$ has at most one neighbor in $X.$ It follows that ${\rm deg}_{G_1}(v)={\rm deg}_{G}(v).$

Denote $H=G_1+x_1y.$ If $H$ is not $2$-connected, then $x_1$ is the only cut-vertex of $H$ and every block of $H$ contains $x_1.$ In this case we let $B$ be the block of $H$ containing $y.$ If $H$ is $2$-connected, we let $B=H.$

If $|B|\ge 3,$ then $B$ is $2$-connected. Now $B-x_1$ is triangle-free and $(B,x_1,y)$ is a $2$-connected rooted graph with minimum degree at least $4$.
By the induction hypothesis, $B$ contains three good $(x_1,y)$-paths, which yield $(u_i,y)$-path $P_i$ in $G$ with $u_i\in X$ for $i=1,2,3$ such that $P_1,P_2$ and $P_3$ are three good paths.
Then $xu_1\cup P_1,$ $xu_2\cup P_2$ and $xu_3\cup P_3$
are three good $(x,y)$-paths in $G$, a contradiction.

Therefore $|B|=2;$ i.e. $B$ is an edge. Then $N_G(y)\subseteq X.$ Since $\deg_G(x)\leq \deg_G(y)$, we conclude that $N_G(y)=X$. Since $G-x$ is triangle-free and $(G,x,y)$ has minimum degree at least $4,$ we have $V(G)\neq X\cup \{y\}.$ So there exists a component $D$ of $G-X$ containing none of $x$ and $y$.

Since $G$ is $2$-connected, we have $|N_G(D)|\geq 2$.
Fix a vertex $u$ in $N_G(D)$.
Let $G_2$ be the graph obtained from $G[N_G[D]]$ by contracting $N_G(D)\setminus\{u\}$ into a new vertex $x_2$.
Since $G-y$ does not have a $4$-cycle containing $x,$ $(G_2,x_2,u)$ is a $2$-connected rooted graph with minimum degree at least $4$. Clearly, $G_2-x_2$ is triangle-free.
By the induction hypothesis, $G_2$ contains three good $(x_2,u)$-paths, which yield $(u_i,u)$-path $P_i$ in $G-\{x,y\}$ with $u_i\in X$ for $i=1,2,3$ such that $P_1,P_2$ and $P_3$ are three good paths. Then
$xu\cup P_1\cup u_1y,$ $xu\cup P_2\cup u_2y$ and $xu\cup P_3\cup u_3y$ are three good $(x,y)$-paths in $G$, a contradiction.
The proof of Claim~3 is complete.

{\bf Claim 4.} There exists a positive integer $s$ and an induced complete bipartite subgraph $Q$ with bipartition $(Q_1,Q_2)$ in $G$ satisfying that
	\begin{enumerate}
		\item $x\in Q_2,\, y\notin V(Q), \, |Q_1|\geq|Q_2|=s+1\geq 2$, and
		\item for every $v\in V(G)\setminus (V(Q)\cup\{y\})$, we have $|N_G(v)\cap Q|\leq s+1$. In particular,
		$|N_G(v)\cap Q_1|\leq s+1$ and $|N_G(v)\cap Q_2|\leq s$.
	\end{enumerate}

By Claim~3, there exists a 4-cycle in $G-y$ containing $x$.
Thus there exists a complete bipartite subgraph $Q$ of $G-y$ with bipartition $(Q_1,Q_2)$ such that $x\in Q_2,$ $y\notin V(Q)$ and $|Q_1|\geq|Q_2|\geq 2$.
We choose $Q$ so that $|Q_2|$ is maximum and subject to this, $|Q_1|$ is maximum.
Let $s$ be the positive integer such that $\lvert Q_2 \rvert = s+1$.

By the choice of $Q$, for every $v\in V(G)\setminus (V(Q)\cup\{y\})$, we have $|N_G(v)\cap Q_1|\leq s+1$ and $|N_G(v)\cap Q_2|\leq s$.
Since $G$ is triangle-free, $Q$ is an induced subgraph in $G$,
and for every $v\in V(G)\setminus (V(Q)\cup\{y\}),$ $v$ cannot be adjacent to both $Q_1$ and $Q_2.$ Hence $|N_G(v)\cap Q|\leq s+1.$
The proof of Claim~4 is complete.

In the remainder of this proof,
$Q$ and $s$ denote the induced complete bipartite subgraph and
the positive integer guaranteed by Claim 4. Additionally, let $C$ be the component of $G-V(Q)$
that contains $y.$

We will now proceed with the proof by dividing it into two cases: $|C|=1$ and $|C|\geq 2$.
In both cases, we will derive a contradiction, thereby showing that $G$ cannot be a counterexample. This will complete the proof of Lemma 5.

{\bf Case 1.} $|C|=1.$

In this case we have $V(C)=\{y\}$.
Recall that by our assumption, $xy \not \in E(G).$ Since $G$ is triangle-free, $y$ is adjacent to exactly one of $Q_1$ and $Q_2$.
Since ${\rm deg}_G(y)\geq {\rm deg}_G(x)$, we deduce that $N_G(x)=N_G(y)=Q_1$ and so $G[V(Q)\cup \{y\}]$ is a complete bipartite graph.
If $s\geq 2$, then $G[V(Q)\cup \{y\}]$ contains three good $(x,y)$-paths of lengths $2,4,6$, respectively, a contradiction. Hence $s=1,$ and we let $\{x_0\}=V(Q_2)\setminus \{x\}.$

We observe that $V(G)\neq V(Q)\cup\{y\}$, for otherwise every vertex in $Q_1$ would have degree $3$ in $G$, a contradiction.
Hence there exists a component in $G-(V(Q)\cup \{y\})$, say $D.$
By Claim~4, we have $|N_G(v)\cap Q|\le 2$ for every $v\in V(D).$ Combining the fact that $\deg_G(v)\ge 4$ for every $v\in V(D)$ with $N_G(y)=Q_1,$ we have $\delta(D)\ge 2.$ This implies that $|V(D)|\geq 3$ and every end-block of $D$ is $2$-connected.
One readily observes that $N_G(D)\cap Q_1\neq\emptyset,$ for otherwise $x_0$ would be a cut-vertex in $G$, contradicting the fact that $G$ is $2$-connected.

{\bf Claim 5.} $|Q_1|\ge 3.$

To the contrary, suppose that $|Q_1|=2.$ Denote by $Q_1=\{u,v\}.$
Since $N_G(x)=N_G(y)=Q_1$, it is easy to check that $(G-\{x,y\},u,v)$ is a $2$-connected rooted graph with minimum degree at least $4$ and $G-\{x,y\}-u$ is triangle-free.
By the induction hypothesis, there are three good $(u,v)$-paths in $G-\{x,y\}$,
which can be easily extended to three good $(x,y)$-paths in $G$, a contradiction.
The proof of Claim~5 is complete.

{\bf Claim 6.} $x_0\in N_G(D).$

Suppose to the contrary that $x_0\notin N_G(D).$
Since $G$ is $2$-connected and $x \notin N_G(D)$, we have $|N_G(D)\cap Q_1|\geq2$.
Let $u_1$ be a vertex in $N_G(D)\cap Q_1$.
Let $G_1$ be the graph obtained from $G[N_G[D]]$ by contracting $N_G(D)\cap (Q_1\setminus \{u_1\})$ into a new vertex $v_1$.
Since $\delta(D)\ge 2,$ the rooted graph $(G_1,v_1,u_1)$ has minimum degree at least $3.$ One readily observes that $(G_1,v_1,u_1)$ is a $2$-connected graph and $G_1-v_1$ is triangle-free. By Lemma 4, $G_1$ contains two nice $(u_1,v_1)$-paths. Hence $G-\{x,y\}$ contains two nice paths $P_i$ from $u_1$ to some vertex $p_i\in V(Q_1)\setminus \{u_1\}$ internally disjoint from $V(Q)$ for $i=1,2.$ Assume that $||P_1||=||P_2||+2.$
By Claim~5, there exists a vertex in $Q_1\setminus \{p_1,u_1\},$ say $x_0'.$
Hence $yu_1P_2p_2x,$ $yu_1P_1p_1x,$ and $yu_1P_1p_1x_0x_0'x,$ are three good $(x,y)$-paths, a contradiction.
The proof of Claim~6 is complete.

{\bf Claim 7.}
There is a matching of size two in $G$ between $V(D)$ and $Q_1$.

Suppose not. Then either $|N_G(D)\cap Q_1|=1$ or $|N_G(Q_1)\cap V(D)|=1$.
In the former case, let $u_2=w_2$ be the unique vertex in $N_G(D)\cap Q_1$;
in the latter case, let $u_2$ be the unique vertex in $N_G(Q_1)\cap V(D)$ and let $w_2$ be a vertex in $Q_1$ adjacent in $G$ to $u_2$.

Recall that $x_0\in N_G(D)$. Let $G_2=G[D\cup \{u_2,x_0\}]$. Then $(G_2,u_2,x_0)$ is a $2$-connected rooted graph with minimum degree at least $4$ and $G_2-u_2$ is triangle-free.
By the induction hypothesis, $G_2$ contains three good $(u_2,x_0)$-paths.
Hence, $G-y$ contains three good $(u_2,x_0)$-paths $P_i$ internally disjoint from $V(Q)$ for $i=1,2,3$. Let $x_0'$ be a vertex in $Q_1\setminus \{w_2\}.$ Concatenating $P_i$ with $yw_2u_2$ and $x_0x_0'x$,
we obtain three good $(x,y)$-paths in $G$, a contradiction.
The proof of Claim~7 is complete.

{\bf Claim 8.} Either $D$ is $2$-connected, or every end-block $B$ of $D$ with cut-vertex $b$ satisfies that $N_G(B-b)\cap Q_1\neq \emptyset.$

To the contrary, suppose that $D$ is not $2$-connected and there exists an end-block $B$ of $D$ with cut-vertex $b$ such that $N_G(B-b)\cap Q_1=\emptyset.$
This implies that $N_G(B-b)\cap V(Q)=\{x_0\},$ as $G$ is $2$-connected.
Recall that every end-block of $D$ is $2$-connected.
So $B$ is $2$-connected.
Let $G_3=G[V(B)\cup \{x_0\}]$.
Then $(G_3,b,x_0)$ is a $2$-connected rooted graph with minimum degree at least $4$ and $G_3-b$ is triangle-free.
By the induction hypothesis, $G_3$ contains three good $(b,x_0)$-paths.
Hence, $G$ contains three good $(b,x_0)$-paths $P_i$ internally disjoint from $V(Q)$ for $i=1,2,3.$

Since $N_G(D)\cap Q_1\neq \emptyset$,
there exists a path $R$ in $G[(D-V(B-b))\cup Q_1]$ from $b$ to some vertex $a\in Q_1$ internally disjoint from $V(B)\cup V(Q)$. Let $x_0'$ be a vertex in $Q_1\setminus \{a\}.$
Concatenating $P_i$ with $R$ and $ay,x_0x_0'x$, we obtain three good $(x,y)$-paths in $G$, a contradiction.
The proof of Claim~8 is complete.

By Claim 7, there exists a matching $M$ of size two in $G$ between $V(D)$ and $Q_1.$
So there exists a vertex $u_4 \in N_G(D)\cap Q_1$ incident with an edge in $M$ such that $N_G(D)\cap (Q_1\setminus \{u_4\})\neq\emptyset$.

Let $G_4$ be the graph obtained from $G[V(D)\cup (N_G(D)\cap Q_1)]$ by contracting $N_G(D)\cap (Q_1\setminus \{u_4\})$ into a new vertex $v_4$.
Since $M$ is a matching of size two in $G$ between $V(D)$ and $Q_1$, if $D$ is $2$-connected, then $(G_4,u_4,v_4)$ is $2$-connected; if $D$ is not $2$-connected, then by Claim 8, every end-block of $D$ has a non-cutvertex adjacent in $G_4$ to one of $u_4,v_4$, so $(G_4,u_4,v_4)$ is $2$-connected.

Recall that $\delta(D)\ge 2$ and $s=1.$ Then $(G_4,v_4,u_4)$ has minimum degree at least three. Clearly, $G_4-v_4$ is triangle-free.
By Lemma 4, there exist two nice $(u_4,v_4)$-paths in $G_4.$
Hence, $G-y$ contains two nice paths $P_i$ from $u_4$ to $p_i\in Q_1\setminus\{u_4\}$ internally disjoint from $V(Q)$ for $i=1,2$. Assume that $||P_1||=||P_2||+2.$
By Claim 5, $|Q_1|\geq 3$ and so we let $x_0'$ be a vertex in $Q_1\setminus \{u_4,p_1\}.$
Then $yu_4P_2p_2x$, $yu_4P_1p_1x$ and $yu_4P_1p_1x_0x_0'x$ are three good $(x,y)$-paths, a contradiction.
This completes the proof of Case 1.

{\bf Case 2.} $|C|\geq 2.$

We first show that $s=1$ or $s=2.$
Suppose to the contrary that $s\ge 3.$ Since $G$ is $2$-connected, $C$ has a neighbor in $Q$, say $u$.
Thus $u\in Q_1\cup Q_2.$  No matter where $u$ lies, we can find three good $(u,x)$-paths in $Q$. Concatenating a fixed  $(u, y)$-path in $G[V(C)\cup \{u\}]$, we obtain three good $(x,y)$-paths in $G,$ a contradiction. So $s=1$ or $s=2.$

{\bf Claim 9.} No vertex in $C-y$ has degree one in $C$.

Suppose to the contrary that there exists $v\in V(C-y)$ with degree one in $C$.
By Claim~4, we have $s+1\geq |N_G(v)\cap V(Q)|\geq 3,$ and so $s=2.$
If $N_G(v)\cap Q_1=\emptyset$, then $N_G(v)\cap V(Q)\subseteq Q_2.$ By Claim 4, $s\geq |N_G(v)\cap V(Q_2)|\geq 3,$ a contradiction.
Hence $N_G(v)\cap Q_1\neq\emptyset.$ Then there are three good $(x,v)$-paths in $G[V(Q)\cup\{v\}]$ of lengths $2,4,6.$
By concatenating each of these path with a $(v, y)$-path in $C$, we obtain three good $(x,y)$-paths in $G$, a contradiction.
The proof of Claim~9 is complete.

By Claim~9, every end-block of $C$ is $2$-connected, except possibly an end-block whose vertex set consists of $y$ and its unique neighbor in $C$.

We say a block $B$ of $C$ is a {\it feasible} block if it is an end-block of $C$ such that either $B=C$ or $y\notin V(B)\setminus\{b\}$ where $b$ is the cut-vertex of $C$
in $B.$  Note that feasible blocks exist, since either $C$ has no cut-vertex, or $C$ contains at least two end-blocks.

Let $B$ be an arbitrary feasible block of $C.$
If $C$ is $2$-connected, then let $b=y$; otherwise let $b$ be the cut-vertex of $C$ contained in $B$.

{\bf Claim 10.} $N_G(B-b)\subseteq Q_2\cup \{b\}$.

Suppose to the contrary that $N_G(B-b)\cap Q_1\neq\emptyset$.
Let $G_1$ be the graph obtained from $G[V(B)\cup (N_G(B-b)\cap Q_1)]$ by contracting $N_G(B-b)\cap Q_1$ into a new vertex $x_1$.
So $(G_1,x_1,b)$ is a $2$-connected rooted graph and $G_1-x_1$ is triangle-free.
Recall that $s=1$ or $s=2.$

Suppose that $s=1.$ By Claim 4, we have $|N_G(v)\cap Q_2|\le 1$ for $v\in V(B-b).$ 
It follows that $(G_1,x_1,b)$ has minimum degree at least three. By Lemma 4, $G_1$ contains two nice $(x_1,b)$-paths. Therefore, there are two nice paths $P_i$ from some vertex $p_i\in N_G(B-b)\cap Q_1$ to $b$ internally disjoint from $V(Q)$ for $i=1,2.$ Also $Q$ contains two $(x,p_i)$-paths of lengths $1,3.$ By concatenating each of these paths with $P_i$ and a fixed $(b,y)$-path in $C-V(B-b),$ we obtain three good $(x,y)$-paths in $G,$ a contradiction.

Suppose that $s=2.$ Let $u$ be a vertex in $N_G(B-b)\cap Q_1$ and let $R$ be a $(u, y)$-path  in $G[V(C)\cup \{u\}].$
One readily observes that $Q$ contains three good $(u,x)$-paths of lengths $1,3,5.$ By concatenating these three paths with $R$, we obtain three good $(x,y)$-paths in $G,$ a contradiction.
The proof of Claim~10 is complete.

{\bf Claim 11.} $s=1$ and $N_G(B-b)\cap V(Q)=Q_2$.

We first show that $|N_G(B-b)\cap V(Q)|\ge 2.$ Suppose not. Let $x'$ be the unique vertex of $N_G(B-b)\cap V(Q).$ Such a vertex $x'$ exists because $G$ is $2$-connected and it is possible that $x'=x.$
Then $(G[N_G[B]],x',b)$ is a $2$-connected rooted graph with minimum degree at least $4$ and $G[N_G[B]]-x'$ is triangle-free, so by the induction hypothesis, $G[V(B)\cup\{x'\}]$ contains three good $(x',b)$-paths.
Let $R$ be a $(b,y)$-path in $C-V(B-b)$ and let $R'$ be a $(x,x')$-path in $Q$ .
Hence concatenating these three good $(x',b)$-paths with $R$ and $R'$ leads to three good $(x,y)$-paths in $G$, a contradiction.
Thus $|N_G(B-b)\cap V(Q)|\ge 2.$ By Claim 10, this implies that $N_G(B-b)\cap (Q_2\setminus\{x\}) \neq \emptyset.$

Suppose that $s=2.$
Let $G_2$ be the graph obtained from $G[V(B)\cup (N_G(B-b)\cap (Q_2\setminus\{x\}))]$ by contracting $N_G(B-b)\cap (Q_2\setminus\{x\})$ into a new vertex $x_2$.
Recall that $C-y$ has no vertex with degree one in $C.$ One readily observes that $(G_2,x_2,b)$ is a $2$-connected rooted graph with minimum degree at least $3$ and $G_2-x_2$ is triangle-free. By Lemma 4, $G_2$ has two nice $(x_2,b)$-paths.
So $G$ contains two nice paths $P_i$ from some vertex $p_i\in N_G(B-b)\cap (Q_2\setminus\{x\})$ to $b$ internally disjoint from $V(Q)$ for $i=1,2$.
Also $Q$ contains two nice $(x,p_i)$-paths of lengths $2,4,$ for $i=1,2.$
By concatenating each of these paths with $P_i$ and $R$, where $R$ is a $(b,y)$-path in $C-V(B-b),$ we obtain three good $(x,y)$-paths, a contradiction.
This shows that $s =1.$ Combining this fact with the inequality $|N_G(B-b)\cap V(Q)|\ge 2,$ we have $N_G(B-b)\cap V(Q)=Q_2$.
The proof of Claim~11 is complete.

Let $a$ be the unique vertex of $Q_2\setminus \{x\}.$

{\bf Subcase 2.1.} $N_G(C-y)\cap Q_1=\emptyset$.

Since $N_G(B-b)\cap V(Q)=\{x,a\}$ and each vertex of $B-b$ has degree at least two in $B,$ $(G[V(B)\cup\{a\}],a,b)$ is a $2$-connected rooted  graph with minimum degree at least three. By Lemma 4, $G[V(B)\cup\{a\}]$ contains two nice $(a,b)$-paths $P_1,P_2.$ Let $Y$ be a $(b, y)$-path  in $C-V(B-b)$.

For any $v\in Q_1$, if $N_G(v)\subseteq Q_2\cup\{y\}$, then the degree of $v$ in $G$ is at most three, a contradiction.
Therefore, there exists a component $D$ of $G-V(Q\cup C)$ adjacent to $v$.

Since $N_G(C-y)\cap Q_1=\emptyset$, we have $N_G(Q_1) \cap V(C) \subseteq \{y\}$.
So $(G-V(C),x,a)$ is a $2$-connected rooted graph with minimum degree at least three and $G-V(C)-x$ is triangle-free.
By Lemma 4, there are two nice $(x,a)$-paths $R_1,R_2$ in $G-V(C)$.
Then $R_i\cup P_j\cup Y$ for all $i,j\in \{1,2\}$ give three good $(x,y)$-paths, a contradiction.

{\bf Subcase 2.2.} $N_G(C-y)\cap Q_1\neq\emptyset$.

If $C$ is $2$-connected, then $C=B$ and $y=b$, which contradicts $N_G(B-b)\cap V(Q)=\{x,a\}$.
So $C$ is not $2$-connected.
Let $B_1,B_2,\ldots, B_t$ be all the end-blocks of $C$ with cut-vertices $b_1,b_2,\ldots,b_t$, respectively.
Clearly, $t\geq 2$.

Suppose that $y\notin\bigcup_{i=1}^t(V(B_i)\setminus\{b_i\})$.
So for every $i \in \{1,2,\ldots,t\},$ $B_i$ is a feasible block, and hence $N_G(B_i-b_i) \cap V(Q)=\{x,a\}$.
Since $N_G(C-y)\cap Q_1\neq\emptyset$,
there is a vertex $w$ in $V(C)\setminus \left(\bigcup_{i=1}^t(V(B_i)\setminus\{b_i\}) \cup \{y\}\right)$ such that $N_G(w)\cap Q_1\neq\emptyset$.
Let $z$ be a vertex in $N_G(w)\cap Q_1$.
Consider the block-cutvertex tree of $C.$ There exist two end-blocks $B_m, B_n$ for $1\leq m<n\leq t$,
such that there are two disjoint paths $L_1, L_2$ from $b_m$ to $w$ and from $b_n$ to $y$ internally disjoint from $V(B_n)\cup V(B_m)$, respectively.

Since $B_m$ and $B_n$ are feasible, $N_G(B_m-b_m) \cap V(Q) = \{x,a\} = N_G(B_n-b_n) \cap V(Q)$.
So both rooted graphs $(G[V(B_m)\cup\{x\}],x,b_m)$ and $(G[V(B_n)\cup\{a\}],a,b_n)$ are $2$-connected and have minimum degree at least $3$. Furthermore, $G[V(B_m)\cup\{x\}]-x$ and $G[V(B_n)\cup\{a\}]-a$ are both triangle-free. By Lemma 4, there are two nice $(x,b_m)$-paths $P_1,P_2$ in $G[V(B_m)\cup\{x\}]$
and two nice $(a,b_n)$-paths $R_1,R_2$ in $G[V(B_n)\cup\{a\}]$.
So the set $\{P_i\cup L_1\cup wza\cup R_j\cup L_2|\, i,j\in \{1,2\}\}$ contains three good $(x,y)$-paths in $G$, a contradiction.

So there exists an end-block, say $B_t$, of $C$ such that $y\in V(B_t)\setminus\{b_t\}$.
We say that a block $H$ of $C$ other than $B_1$ is a {\it hub} if $H$ is $2$-connected and contains at most two cut-vertices of $C$, and every path in $C$ from $b_1$ to $b_t$ contains all cut-vertices of $C$ contained in $V(H)$.

Suppose there exists a hub $B^*$ of $C$.
So there exists a cut-vertex $x^*$ of $C$ contained in $B^*$ such that every path in $C$ from $b_1$ to $V(B^*)$ contains $x^*$.
If $B^*=B_t$, then let $y^*=y$; otherwise, let $y^*$ be the cut-vertex of $C$ contained in $B^*$ such that every path in $C$ from $b_t$ to $V(B^*)$ contains $y^*$.
Let $Z_0$ be a $(b_1,x^*)$-path in $C-(V(B_1-b_1)\cup V(B^*-x^*)),$ and let $Z_1$ be a $(y^*,y)$-path in $C.$
Since $(G[B_1\cup \{x\}],x,b_1)$ is a $2$-connected rooted graph with minimum degree at least three and $G[B_1\cup \{x\}]-x$ is triangle-free, by Lemma 4, $G[B_1\cup \{x\}]$ contains two nice $(x,b_1)$-paths $P_1,P_2.$

If every vertex in $V(B^*)\setminus\{x^*,y^*\}$ has at most one neighbor in $Q$, then $(B^*,x^*,y^*)$ is a $2$-connected rooted graph with minimum degree at least three and $B^*-x^*$ is triangle-free. By Lemma 4, $B^*$ contains two nice $(x^*,y^*)$-paths $R_1,R_2.$ Hence, the set $\{P_i\cup Z_0\cup R_j \cup Z_1|\,  i,j\in \{1, 2\}\}$ contains three good $(x,y)$-paths in $G$, a contradiction.

Therefore some vertex $w\in V(B^*)\setminus\{x^*,y^*\}$ satisfies $|N_G(w)\cap V(Q)|\geq 2$.
Since $s=1$, by Claim~4, we have $|N_G(w)\cap V(Q)|=2$.
Let $u,v$ be the vertices in $N_G(w)\cap V(Q)$.
By Claim~4 and the fact that $G$ is triangle-free, we have $\{u,v\} \subseteq Q_1$.
Hence there are two nice $(x,a)$-paths $L_1=xua$ and $L_2=xuwva.$ Since $(G[B_1\cup \{a\}],a,b_1)$ is a $2$-connected rooted graph with minimum degree at least three and $G[B_1\cup \{a\}]-a$ is triangle-free, by Lemma~4, there exist two nice $(a,b_1)$-paths $N_1,N_2$ in $G[B_1\cup \{a\}]$. Since $B^*$ is $2$-connected, there exists a $(x^*,y^*)$-path $L'$ in $B-w$. Therefore, the set $\{L_i\cup N_j\cup Z_0\cup L' \cup Z_1|\,  i,j\in \{1, 2\}\}$ contains three good $(x,y)$-paths in $G$, a contradiction.

So there exists no hub.
In particular, $B_t$ is not $2$-connected, for otherwise $B_t$ is a hub.
Therefore $B_t=yb_t$ is an edge. So $B_1,...,B_{t-1}$ are the only feasible blocks in $C$.
Recall that $N_G(B_i-b_i)\cap V(Q)=\{a,x\}$ for all $i\in \{1,2,\ldots, t-1\},$ which implies $\deg_G(x)\geq |Q_1|+t-1$.

Since $G$ is triangle-free, we have $\deg_G(y)\leq |Q_1|+1$.
Hence $|Q_1|+t-1\leq \deg_G(x)\leq \deg_G(y)\leq |Q_1|+1$.
That is, $t \leq 2$.
As $t\geq 2$, this forces $t=2$, $\deg_G(x)=\deg_G(y)=|Q_1|+1$.
In other words, there is exactly one end-block $B_1$ of $C$ other than $B_2=yb_2$, $N_G(y)=Q_1\cup \{b_2\}$ and $N_G(x)\subseteq Q_1\cup V(B_1-b_1)$.

Note that the block-cutvertex tree of $C$ is a path.
Since there exists no hub, every block of $C$ other than $B_1$ is an edge.
If $b_1=b_2$, then since $N_G(C-y)\cap Q_1\neq\emptyset$ and $N_G(B_1-b_1)\cap Q_1=\emptyset$, $b_2$ must have a neighbor in $Q_1$.
If $b_1\neq b_2$, then $|N_G(b_2)\cap V(C)|=2$, and since ${\rm deg}_G(b_2)\geq 4$, we have $|N_G(b_2)\cap V(Q)|\geq 2$.
Recall that $N_G(x)\subseteq Q_1\cup V(B_1-b_1),$ so $xb_2\notin E(G)$.
Thus in either case, $b_2$ must have a neighbor $w^*$ in $Q_1$.
But $G[\{y,b_2,w^*\}]$ is a triangle, a contradiction.
This completes the proof of Lemma 5. \hfill $\Box$

The following concept is crucial in our approach.
We say that a cycle $C$ in a connected graph $G$ is {\it non-separating} if $G-V(C)$ is connected.
The proof of the following lemma can be found in [2, p.14, proof of Theorem 2], though it was not formally stated.

{\bf Lemma 6} (Bondy-Vince [2, p.14]). {\it Every $3$-connected nonbipartite graph contains a non-separating induced odd cycle.}

We also need the following lemma on non-separating odd cycles due to Liu and Ma [9] which is a slight modification of a result of Fan [4].

{\bf Lemma 7} (Liu-Ma [9, Lemma 5.1, p.88]). {\it Let $G$ be a graph with minimum degree at least four.
If $G$ contains a non-separating induced odd cycle, then $G$ contains a non-separating induced odd cycle $C$, denoted by $v_0v_1...v_{2s}v_0$,  such that either

(1) $C$ is a triangle, or

(2) for every non-cut-vertex $v$ of $G-V(C)$, $\lvert N_G(v)\cap V(C) \rvert \le 2$, and equality holds if and only if $N_G(v)\cap V(C)=\{v_i,v_{i+2}\}$ for some $i$,
where all subscripts of $v_i$ are read modulo $2s+1.$
}

When the graph contains a triangle, the conclusion of Theorem 3 has already been proved by Gao, Huo and Ma as follows.

{\bf Lemma 8} (Gao-Huo-Ma [6, Theorem 4.1, p.2320]). {\it Let $k \ge 2$ be an integer and let $G$ be a $2$-connected graph containing a triangle. If $G$ has minimum degree
at least $k$ and order at least $k+2$, then $G$ contains $k$ cycles of consecutive lengths.}

Now we are ready to prove the main result.

{\bf Proof of Theorem 3.} The case $k\ge 6$ of Theorem 3 has been proved in [8, Theorem 1.2]. We need only to consider the cases $k=4$ and $k=5.$

By Lemma~8, it suffices to consider the case when $G$ is triangle-free and we make this assumption. Combining Lemmas 6 and 7, we
deduce that $G$ contains a non-separating induced odd cycle $C=v_0v_1...v_{2s}v_0$ that satisfies the property (2) of Lemma 7. Let $G_1=G-V(C).$

The following claim can be found in [8, Section 3].

{\bf Claim 1.} Given $k \ge 4$,  let $B =G_1$ if  $G_1$ is $2$‐connected; otherwise, let $B$ be an  end-block of $G_1$ with cut-vertex $x$. If there exists $v\in V(B)\setminus\{x\}$ such that $|N_G(v) \cap V(C)|=2$, then $G$  contains $k$ cycles of consecutive lengths.

Depending on whether  $G_1$ is $2$-connected,  we divide the rest of the proof into two cases.

{\bf Case 1. } $G_1$  is $2$-connected.

By Claim~1,
 we can assume that every vertex $v\in V(G_1)$ satisfies $|N_G(v) \cap V(C)| \le 1$ which implies that the minimum degree of $G_1$ is at least $k-1$. We choose a vertex $v\in N_G(v_0)\cap V(G_1)$ and a vertex $u\in N_G(v_{s})\cap V(G_1)$ such that $u\neq v$. Note that $vv_0v_1v_2\dots v_{s}u$ and $vv_{0}v_{2s}v_{2s-1}\dots v_{s}u$ are two
 $(v,u)$-paths of lengths $s+2$ and $s+3$, denoted by $Q_1$  and $Q_2$, respectively.
 Applying Lemmas 4 and 5 to $G_1$ with $u,v\in V(G_1)$, we deduce that $G_1$ contains $k-2$ nice or good $(u,v)$-paths $P_1, P_2, \dots, P_{k-2}$
 according as $k=4$ or $k=5.$  By concatenating the paths $Q_i$ and $P_j$ we obtain $2k-4$ cycles $Q_i\cup P_j$ for $i=1,2$ and $j=1,\ldots, k-2$ where $k=4$ or $k=5.$
 Clearly, among these $2k-4$ cycles, there exist  $k$ cycles of consecutive lengths.

{\bf Case 2. } $G_1$ is not $2$-connected.

By Claim~1, we can assume that for any end-block $B$ of $G_1$, every vertex of $B$ other than the cut-vertex has degree at least $k-1$.
Let $D_1$ be an end-block with cut-vertex $x$ of $G_1$ and $G_2=G_1-(V(D_1)\setminus \{x\})$.
Note that both $D_1$ and $G_2$ have orders at least $3.$
The following claim can be found in [8, Section 3].

{\bf Claim 2.} There exists an integer $i$ such that $N_G(v_i)\cap (V(D_1)\setminus \{x\})\neq \emptyset $ and $N_G(v_{i+s})\cap V(G_2)\neq \emptyset.$

Now let $i'$ be an integer such that $N_G(v_{i'})\cap (V(D_1)\setminus \{x\})$ contains a vertex $v$ and $N_G(v_{i'+s})\cap V(G_2)$ contains a vertex $u.$
Note that $D_1$ is $2$-connected and every vertex of $D_1$ other than $x$ has degree at least $k-1$. By Lemmas 4 and 5, $D_1$ contains $k-2$ nice or good $(x,v)$-paths $P_1, \dots, P_{k-2}.$

Let $Q_1=vv_{i'}v_{i'+1}\dots v_{i'+s}u$ and $Q_2=uv_{i'+s}v_{i'+s+1}\dots v_{i'}v.$
Then $Q_1$ and $Q_2$ are two $(v,u)$-paths with consecutive lengths in $G[V(C)\cup \{u,v\}].$
Let $T$ be a $(u,x)$-path in $G_2$, possibly $x=u$.
By concatenating the paths $Q_i,$ $T$ and $P_j$ we obtain $2k-4$ cycles $Q_i\cup T\cup P_j$ for $i=1,2$ and $j=1,\ldots, k-2$ where $k=4$ or $k=5.$
 Clearly, among these $2k-4$ cycles, there exist  $k$ cycles of consecutive lengths.
This completes the proof of Theorem 3. \hfill $\Box$


\vskip 5mm
{\bf Acknowledgement.} This research  was supported by the NSFC grant 12271170 and Science and Technology Commission of Shanghai Municipality
 grant 22DZ2229014.

\end{document}